\documentclass [10pt,leqno,a4paper]{article}
\usepackage{amssymb}
\usepackage{enumerate}
\usepackage{amsmath}
\usepackage{amsfonts,mathrsfs}
\usepackage{amsthm}
\usepackage{epsfig}

\newcounter{teoremaganso}

\renewenvironment{abstract}{\small\quotation\noindent
 {\bfseries \abstractname .}}{\endquotation \par}

\newenvironment{trivlistparm}[2]{\trivlist
\item[{#1 #2}]}{\endtrivlist}
\newenvironment{proclama}[1]{\trivlistparm{\bfseries}{#1}\itshape}{\endtrivlistparm}
\newenvironment{prooftext}[1]{\trivlistparm{\bfseries}{#1}}{\Qed\endtrivlistparm}
\newenvironment{prova}{\trivlistparm{\bfseries}{Proof.}}{\Qed\endtrivlistparm}

\catcode`\@=11

\def\resetthefootnote{\renewcommand{\thefootnote}{\@arabic\c@footnote} }
\def\@principiremex#1{\trivlist
 \item[\hskip \labelsep{\bfseries #1\ \thetheo}]\ignorespaces}
\def\opar@principiremex#1[#2]{\trivlist
 \item[\hskip \labelsep{\bfseries #1\ \thetheo\ (#2)}]\ignorespaces}

\newcommand{\newTHEOremrom}[2]{\newenvironment{#1}{\refstepcounter{theo}\@ifnextchar[{\opar@principiremex{#2}}{\@principiremex{#2}}
  }{\qedB\endtrivlist}}
\catcode`\@=12 \DeclareMathSymbol{\square}{\mathord}{AMSa}{"03}
\newcommand{\qedB}{\nopagebreak\hspace*{\fill}$\square$\par}
\newcommand{\Qed}{\nopagebreak\hspace*{\fill}{\vrule width6pt height6pt depth0pt}\par}

\newtheorem {theo} {Theorem} [section]
\newtheorem {prop} [theo] {Proposition}

\newtheorem {lem} [theo] {Lemma}
\newtheorem {bigtheo} [teoremaganso] {Theorem}

\newTHEOremrom {defi} {Definition}
\newTHEOremrom {obs} {Remark}
\newTHEOremrom {ex} {Example}

\newcommand{\refc}[1]{\mbox{$(\ref{#1})$}}
\newcommand{\secc}[1]{Section~\ref{#1}}
\newcommand{\teoc}[1]{Theorem~\ref{#1}}
\newcommand{\propc}[1]{Proposition~\ref{#1}}

\newcommand{\lemc}[1]{Lemma~\ref{#1}}
\newcommand{\defic}[1]{Definition~\ref{#1}}
\newcommand{\obsc}[1]{Remark~\ref{#1}}
\newcommand{\exc}[1]{Example~\ref{#1}}
\newcommand{\figc}[1]{Figure~\ref{#1}}

\newcommand{\N}{\ensuremath{\mathbb{N}}}

\newcommand{\R}{\ensuremath{\mathbb{R}}}

\newcommand{\C}{\ensuremath{\mathbb{C}}}
\newcommand{\ii}{\ensuremath{(x_{\ell},x_r)}}

\newcommand{\op}{\ensuremath{\mbox{\rm o}}}

\def\map#1#2#3{\mbox{${#1}\!:{#2}\longrightarrow{#3}$}}
\newcommand{\dsp}{\displaystyle}
\newcommand{\dspt}[1]{\displaystyle\textstyle{#1}}

\begin{document}

\title{\textbf{A Chebyshev criterion for Abelian integrals}
\footnotetext{The first author is partially
supported by the MEC/FEDER grant MTM2005-06098-C02-02. The second
author by the MEC/FEDER grants MTM2005-02139 and MTM2005-06098 and
the CIRIT grant 2005SGR-00550. The third author by the MEC/FEDER
grant MTM2005-06098-C02-01 and the CIRIT grant 2005SGR-00550.}}

\author{M. Grau, F. Ma\~{n}osas and J. Villadelprat
\\*[.1truecm]
{\small \textsl{Departament de Matem\`{a}tica, }}
\\*[-.05truecm]
{\small \textsl{Universitat de Lleida, Lleida, Spain}}
\\*[.1truecm]
{\small \textsl{Departament de Matem\`{a}tiques}}
\\*[-.05truecm]
{\small \textsl{Universitat Aut\`{o}noma de Barcelona, Barcelona, Spain}}
\\*[.1truecm]
{\small \textsl{Departament d'Enginyeria Inform\`{a}tica i Matem\`{a}tiques,}}
\\*[-.05truecm]
{\small \textsl{Universitat Rovira i Virgili, Tarragona, Spain}}}

\date{}

\maketitle

\begin{abstract}
We present a criterion that provides an easy sufficient condition in
order that a collection of Abelian integrals has the Chebyshev
property. This condition involves the functions in the integrand of
the Abelian integrals and can be checked, in many cases, in a purely
algebraic way. By using this criterion, several known results are
obtained in a shorter way and some new results, which could not be
tackled by the known standard methods, can also be deduced.
\end{abstract}

{\small{\noindent }

\section{Introduction and statement of the result}\label{Sec1}

The second part of \emph{Hilbert's 16th problem} \cite{Hil} asks
about the maximum number and location of limit cycles of a planar
polynomial vector fields of degree $d.$ Solving this problem, even
in the case $d=2,$ seems to be out of reach at the present state of
knowledge (see the works of Ilyashenko \cite{Ily} and Li Jibin
\cite{Jibin} for a survey of the recent results on the subject). Our
paper is concerned with a weaker version of this problem, the
so-called \emph{infinitesimal Hilbert's 16th problem}, proposed by
Arnold \cite{Arnoldp}. Let $\omega$ be a real 1-form with polynomial
coefficients of degree at most $d.$ Consider a real polynomial $H$
of degree $d+1$ in the plane. A closed connected component of a
level curve $H=h$ is denoted by~$\gamma_h$ and called an \emph{oval}
of~$H.$ These ovals form continuous families (see \figc{dib2}) and
the infinitesimal Hilbert's 16th problem is to find an upper bound
$V(d)$ of the number of real zeros of the \emph{Abelian integral}
 \begin{equation}\label{Abeliana}
 I(h)=\int_{\gamma_h}\omega.
 \end{equation}
The bound should be uniform with respect to the choice of the
polynomial~$H,$ the family of ovals~$\{\gamma_h\}$ and the
form~$\omega.$ It should depend on the degree $d$ only. (In the
literature an Abelian integral is usually the integral of a
rational 1-form over a continuous family of algebraic ovals.
Throughout the paper, by an abuse of language, we use the name
Abelian integral also in case the functions are analytic.)

Zeros of Abelian integrals are related to limit cycles in the
following way. Consider a small deformation of a Hamiltonian
vector field $X_{\varepsilon}=X_H+\varepsilon Y,$ where
 \[
  X_H=-H_y\partial_x+H_x\partial_x\,\mbox{ and }\,Y=P\partial_x+Q\partial_y.
 \]
Then, see \cite{Ily,Jibin} for details, the first approximation in
$\varepsilon$ of the displacement function of the Poincar\'{e} map of
$X_{\varepsilon}$ is given by~\refc{Abeliana} with $\omega=Pdy-Qdx.$
Hence the number of isolated zeros of $I(h),$ counted with
multiplicities, provides an upper bound for the number of ovals of
$H$ that generate limit cycles of $X_{\varepsilon}$ for
$\varepsilon\approx 0.$ The coefficients of~$P$ and~$Q$ are
considered as parameters of the problem and so the function $I(h)$
splits as a linear combination
\[
 \alpha_0I_0(h)+\alpha_1I_1(h)+\ldots+\alpha_{n-1}I_{n-1}(h),
\]
where $\alpha_k$ depends on the initial parameters and $I_k(h)$ is
an Abelian integral with either $\omega=x^iy^jdx$ or
$\omega=x^iy^jdy$. (In fact it is easy to see, using integration by
parts, that only one type of these 1-forms needs to be considered.)
Therefore the problem is equivalent to find an upper bound for the
number of isolated zeros of any function belonging to the vector
space generated by $I_k(h)$ for $k=0,1,\ldots,n-1.$ This problem is
strongly related to showing that the basis of the previous vector
space is a Chebyshev system. In fact, the great majority of papers
studying concrete problems on the subject show this kind of
property.

In this paper we focus on the case in which $H$ has separated
variables, i.e., $H(x,y)=\Phi(x)+\Psi(y),$ and as a byproduct we
obtain a result for the case $H(x,y)=A(x)+B(x)y^{2m}$ as well. We
suppose in addition that
\begin{equation*}
  I_i(h)=\int_{\gamma_h}f_i(x)g(y)dx,\ \mbox{ for $i=0,1,\ldots,n-1,$}
 \end{equation*}
where $f_0,f_1,\ldots,f_{n-1}$ and $g$ are \emph{analytic
functions.} (Note that the function depending on $y$ is the same for
all the 1-forms. In the problems studied in the literature, the
original family of Abelian integrals can be usually reduced to a
family as above.) We will show that, in this case, some Chebyshev
properties on $f_i$ and $g$ (to be specified later on) transfer to
$I_i$ after the integration over the ovals. To fix notation, $H$ is
an analytic function in some open subset of the plane that has a
local minimum at the origin. Then there exists a punctured
neighbourhood $\mathcal P$ of the origin foliated by ovals
$\gamma_h\subset\{H(x,y)=h\}.$ We fix that $H(0,0)=0$ and then the
set of ovals~$\gamma_h$ inside this, let us say, \emph{period
annulus,} can be parameterized by the energy levels $h\in (0,h_0)$
for some $h_0\in (0,+\infty]$. In what follows, we shall denote the
projection of $\mathcal P$ on the $x$-axis by $\ii$. Similarly,
$(y_{\ell},y_r)$ is the projection of $\mathcal P$ on the $y$-axis.

\teoc{desacoplat} is our main result and it applies in case that
$H(x,y)=\Phi(x)+\Psi(y).$ It is easy to verify that, under the above
assumptions, $x\Phi'(x)>0$ for any $x\in\ii\setminus\{0\}$ and
$y\Psi'(y)>0$ for any $y\in(y_{\ell},y_r)\setminus \{0\}.$
Then~$\Phi$ and~$\Psi$ must have even multiplicity at $0.$ Thus,
there exist two analytic involutions~$\sigma_1$ and~$\sigma_2$ such
that
\begin{align*}
 &\mbox{$\Phi(x)=\Phi\bigl(\sigma_1(x)\bigr)$ for all $x\in\ii$}
 \intertext{and}
 &\mbox{$\Psi(y)=\Psi\bigl(\sigma_2(y)\bigr)$ for all $y\in (y_{\ell},y_r).$}
\end{align*}
Recall that a mapping $\sigma$ is an \emph{involution} if
$\sigma\circ\sigma=Id$ and $\sigma\neq Id.$ Note that an involution
is a diffeomorphism with a unique fixed point. In our situation we
have that $\sigma_i(0)=0.$ In what follows, given a function
$\kappa,$ we define its \emph{balance} with respect to $\sigma$ as
\[
\mathscr B_{\sigma} \! \bigl(\kappa\bigr)(x) =
\kappa(x)-\kappa\bigl(\sigma(x)\bigr).
\]
For example, if $\sigma=-Id$, then the balance of a function is
twice its odd part.

In the statement of \teoc{desacoplat}, $m$ is related with the
multiplicity of $\Psi$ at $y=0.$ More concretely, we suppose that
$\Psi(y)=ey^{2m}+\op(y^{2m})$ with $e>0.$ In addition, ECT-system
stands for \emph{extended complete Chebyshev} system in the sense
of Marde\v si\'c~\cite{Mar}, see \defic{difi1} for details.

\begin{bigtheo}\label{desacoplat}
Let us consider the Abelian integrals
\begin{equation*}
I_i(h)\, =\, \int_{\gamma_h}f_i(x)g(y)dx,\ \mbox{
$i=0,1,\ldots,n-1,$}
\end{equation*}
where, for each $h\in (0,h_0),$ $\gamma_h$ is the oval surrounding the origin inside the
level curve $\{\Phi(x)+\Psi(y)=h\}.$ Let $\sigma_1$ and $\sigma_2$ be the involutions
associated to $\Phi$ and $\Psi$,  respectively. Setting $g_0=g,$ we define
$g_{i+1}=\frac{g'_i}{\Psi'}$. Then $(I_0,I_1,\ldots,I_{n-1})$ is an ECT-system on
$(0,h_0)$ if the following hypothesis are satisfied:
\begin{enumerate}[$(a)$]
\item $\Bigl(\mathscr B_{\sigma_1}\!\bigl(\frac{f_0}{\Phi'}\bigr),
        \mathscr B_{\sigma_1}\!\bigl(\frac{f_1}{\Phi'}\bigr),\ldots,
        \mathscr B_{\sigma_1}\!\bigl(\frac{f_{n-1}}{\Phi'}\bigr)\Bigr)$ is a CT-system on
        $(0,x_r),$ and
\item $\Bigl(\mathscr B_{\sigma_2}\!(g_0),\mathscr B_{\sigma_2}\!(g_1),\ldots,
        \mathscr B_{\sigma_2}\!(g_{n-1})\Bigr)$ is a CT-system on $(0,y_r)$ and
        $\mathscr B_{\sigma_2}\!(g_0)(y)=\op(y^{2m(n-2)}).$
\end{enumerate}
\end{bigtheo}

To prove the result it is necessary to compute the derivative of
each Abelian integral until order $n-1$. The condition on
$\mathscr B_{\sigma_2}\!(g_0)(y)$ at $y=0$ ensures that the
integral expression of this derivative is convergent, although it
may be improper (see \obsc{impropia}). Let us also point out that,
since $\sigma_2(y)=-y+\op(y),$ this condition is equivalent to
require that $g(y)-g(-y)=\op(y^{2m(n-2)}).$

Our second result deals with those Abelian integrals such that
\[
H(x,y)=A(x)+B(x)y^{2m}\ \mbox{ and $g(y)=y^{2s-1}$ with $s\in\N.$}
\]
Since $H$ has a local minimum at the origin by assumption, $B(0)>0$
and $A$ has a local minimum at $x=0.$ Thus, as before, there exists
an involution $\sigma$ satisfying $A(x)=A\bigl(\sigma(x)\bigr)$ for
all $x\in\ii$.

\begin{bigtheo}\label{quadratic}
Let us consider the Abelian integrals
\[
I_i(h)=\int_{\gamma_h}f_i(x){y^{2s-1}}dx,\ \mbox{
$i=0,1,\ldots,n-1,$}
\]
where, for each $h\in (0,h_0),$ $\gamma_h$ is the oval surrounding the origin inside the
level curve $\{A(x)+B(x)y^{2m}=h\}.$ Let~$\sigma$ be the involution associated to $A$ and
we define
 \[
  \dspt{\ell_i=\mathscr B_{\sigma}\!\left(\frac{f_i}{A'B^{\frac{2s-1}{2m}}}\right).}
 \]
Then $(I_0,I_1,\ldots,I_{n-1})$ is an ECT-system on $(0,h_0)$ if $s>m(n-2)$ and
$\bigl(\ell_0,\ell_1,\ldots,\ell_{n-1}\bigr)$ is a CT-system on $(0,x_r).$
\end{bigtheo}

It is worth noting that although the condition $s>m(n-2)$ is not
fulfilled in some situations, it is possible to obtain a new
Abelian integral for which the corresponding $s$ is large enough
to verify the inequality. The procedure to obtain this new Abelian
integral follows from the application of \lemc{puja}. We refer the
reader to \exc{ex1} in which we explain in detail how to apply
\lemc{puja} to get a new Abelian integral with $s>m(n-2)$.

The applicability of our criteria comes from the fact that the hypothesis requiring some
functions to be a CT-system can be verified by computing Wronskians (see \lemc{lem1}).
This simplifies a lot the problem of showing that a given collection of Abelian integrals
has the Chebyshev property and in some cases it enables to reformulate the problem in a
purely algebraic way (cf. \secc{Sec4}).

In the literature there are a lot of papers dealing with zeros of
Abelian integrals (see for instance
\cite{DR,Han,Gasull,Girard,Gavrilov3,Petrov,Zhao} and references
there in). In many cases, it is essential to show that a collection
of Abelian integral has some kind of Chebyshev property. The
techniques and arguments to tackle these problems are usually very
long and highly non-trivial. For instance, in some papers (e.g.
\cite{DLZ,Horo,Peng}) the authors study the geometrical properties
of the so-called \emph{centroid curve} using that it verifies a
Riccati equation (which is itself deduced from a Picard-Fuchs
system). In other papers (e.g. \cite{Iliev3,Gavrilov1,Gavrilov2}),
the authors use complex analysis and algebraic topology (analytic
continuation, argument principle, monodromy, Picard-Lefschetz
formula, \ldots). Certainly, the criterion that we present here can
not be applied to all the situations (since the Abelian integrals
need to have a specific structure) and, even in case that it is
possible to apply it, sometimes the sufficient condition that we
provide is not verified. However we want to stress that, when it
works, it enables to extremely simplify the solution. To illustrate
this fact, in \secc{Sec4} we reprove with our criterion the main
results of three different papers. We are also convinced that this
criterion will be useful to obtain new results on the issue. In this
direction we tackle the program posed by Gautier, Gavrilov and Iliev
\cite{Iliev3} and we prove their conjecture in four new cases (see
Subsection \ref{Gautier}).

In several papers dealing with zeros of Abelian integrals (see
\cite{DumLi1,DumLi2,DLZ,Peng} for instance), it is applied a
criterion of Li and Zhang~\cite{LiZhang}. This criterion provides a
sufficient condition for the monotonicity of the ratio of two
Abelian integrals. In page~360 of the book of Arnold's problems
\cite{Arnoldp}, the criterion given in \cite{LiZhang} is quoted as a
useful tool that ``despite its seemingly artificial form, it proves
to be working in many independently arising particular cases''. The
translation of the result in~\cite{LiZhang} to the language of
Chebyshev systems and Wronskians shows that it corresponds precisely
to the case $n=2$ of our criteria. Accordingly, using our
formulation, their result becomes very natural: it shows that the
Chebyshev properties of the functions in the 1-form are preserved
after integration. In addition, as a generalization of their result,
we hope that our criteria will be useful in many cases as well.
Finally we remark that, although we suppose that the functions that
we deal with are analytic, our results hold true for smooth
functions with minor changes.

The paper is organized as follows. \secc{Sec2} is devoted to introduce the definitions
and the notation  that we shall use. In particular we define the different types of
Chebyshev property that we shall deal with and we establish their equivalences with the
continuous and discrete Wronskians (see \lemc{lem1}). Theorems~\ref{desacoplat}
and~\ref{quadratic} are proved in \secc{Sec3}. The main ingredient in the proof of
\teoc{desacoplat} is \propc{prop3}, that provides an integral expression for the
Wronskian of a collection of Abelian integrals. \teoc{quadratic} follows as a corollary
of \teoc{desacoplat}. \secc{Sec4} is devoted to illustrate the application of our
criteria. To this end, in Examples~\ref{ex1}, \ref{ex2} and~\ref{ex3} we reprove the
results of Iliev and Perko~\cite{Iliev3}, Zhao, Liang and Lu~\cite{Zhao} and
Peng~\cite{Peng}, respectively. Apart from showing the simplicity in the application of
the criteria, our aim with these examples is twofold. First, to show that it is not
necessary to know explicitly the involutions that appear in the statements. Second, to
show that it is possible to reformulate the problem in such a way it suffices to check
that some \emph{polynomials} do not vanish. In \secc{Sec4} we also present some new
results concerning the program of Gautier, Gavrilov and Iliev~\cite{Iliev3}. Finally in
the Appendix we give some details about the tools that are used in \secc{Sec4}, namely,
the notion of resultant between two polynomials and Sturm's Theorem.

\section{Chebyshev systems}\label{Sec2}

\begin{defi}\label{difi1}
Let $f_0,f_1,\ldots,f_{n-1}$ be analytic functions on an open interval $L$ of $\R.$
\begin{enumerate}[$(a)$]
\item $(f_0,f_1,\ldots,f_{n-1})$ is a \emph{Chebyshev system} $($in short, T-system$)$ on
      $L$ if any nontrivial linear combination
      \[
       \alpha_0f_0(x)+\alpha_1f_1(x)+\ldots+\alpha_{n-1}f_{n-1}(x)
      \]
      has at most $n-1$ isolated zeros on $L.$
\item $(f_0,f_1,\ldots,f_{n-1})$ is a \emph{complete Chebyshev
      system} $($in short, CT-system$)$ on $L$ if
      $(f_0,f_1,\ldots,f_{k-1})$ is a T-system for all $k=1,2,\ldots,n.$
\item $(f_0,f_1,\ldots,f_{n-1})$ is an \emph{extended complete
      Chebyshev system} $($in short, ECT-system$)$ on $L$ if, for all $k=1,2,\ldots,n,$ any
      nontrivial linear combination
      \[
       \alpha_0f_0(x)+\alpha_1f_1(x)+\ldots+\alpha_{k-1}f_{k-1}(x)
      \]
      has at most $k-1$ isolated zeros on $L$ counted with multiplicities.
\end{enumerate}
(Let us mention that, in these abbreviations, ``T'' stands for
Tchebycheff, which in some sources is the transcription of the
Russian name Chebyshev.)
\end{defi}

It is clear that if $(f_0,f_1,\ldots,f_{n-1})$ is an ECT-system on
$L$, then $(f_0,f_1,\ldots,f_{n-1})$ is a CT-system on $L$.
However, the reverse implication is not true.

\begin{defi}
Let $f_0,f_1,\ldots,f_{k-1}$ be analytic functions on an open
interval $L$ of $\R.$ The \emph{continuous Wronskian} of
$(f_0,f_1,\ldots,f_{k-1})$ at $x\in L$ is
 \[
 W\bigl[f_0,f_1,\cdots,f_{k-1}\bigr](x)=\det\left(f_j^{(i)}(x)\right)_{0\leqslant i,j\leqslant k-1}
   =\left|\begin{array}{ccc}
    f_0(x) & \cdots & f_{k-1} (x) \\
    f'_0(x) & \cdots & f'_{k-1} (x) \\
    & \vdots & \\
    f^{(k-1)}_{0}(x) & \cdots & f^{(k-1)}_{k-1} (x) \\
    \end{array} \right|
  \]
The \emph{discrete Wronskian} of $(f_0,f_1,\ldots,f_{k-1})$ at
$(x_0,x_1,\ldots, x_{k-1})\in L^k$ is
 \[
  D\bigl[f_0,f_1,\cdots,f_{k-1}\bigr](x_0,x_1,\ldots,x_{k-1})
    =\det\bigl(f_j(x_i)\bigr)_{0\leqslant i,j\leqslant k-1}
    =\left|\begin{array}{ccc}
    f_0(x_0) & \cdots & f_{k-1}(x_0) \\
    f_0(x_1) & \cdots & f_{k-1}(x_1) \\
    & \vdots & \\
    f_{0}(x_{k-1}) & \cdots & f_{k-1}(x_{k-1}) \\
    \end{array} \right|
 \]
\end{defi}

For the sake of shortness, given any ``letter'' $x$ and $k\in\N$ we
use the notation
 \[
  x_0,x_1,\ldots,x_{k-1}=\mathbf{x_k}.
 \]
Accordingly, we write
 \begin{align*}
 &W\bigl[f_0,f_1,\cdots,f_{k-1}\bigr](x)=W\bigl[\mathbf{f_k}\bigr](x)\\
 \intertext{and}
 &D\bigl[f_0,f_1,\cdots,f_{k-1}\bigr](x_0,x_1,\ldots,x_{k-1})=D\bigl[\mathbf{f_k}\bigr](\mathbf{x_k})
 \end{align*}
for the continuous and discrete Wronskian, respectively. The
following result is well known (see \cite{Karlin,Mar} for
instance).

\begin{lem}\label{lem1}
The following equivalences hold:
\begin{enumerate}[$(a)$]
 \item $(f_0,f_1,\ldots,f_{n-1})$ is a CT-system on $L$ if, and only if,
       for each $k=1,2,\ldots,n,$
       \[
        D\bigl[\mathbf{f_k}\bigr](\mathbf{x_k})\neq 0\,
        \mbox{ for all $\mathbf{x_k}\in L^k$ such that $x_i\neq x_j$ for $i\neq j.$}
       \]
 \item $(f_0,f_1,\ldots,f_{n-1})$ is an ECT-system on $L$ if, and only if,
       for each $k=1,2,\ldots,n,$
       \[
        W\bigl[\mathbf{f_k}\bigr](x)\neq 0\mbox{ for all $x\in L.$}
       \]
\end{enumerate}
\end{lem}

\section{Proof of the main results}\label{Sec3}

The first part of this section is devoted to prove
\teoc{desacoplat}. Thus, unless we explicitly say the contrary, we
suppose that $H(x,y)=\Phi(x)+\Psi(y)$, where
$\Psi(y)=ey^{2m}+\op(y^{2m})$ with $e>0$, as mentioned before.
Then, there exists a diffeomorphism $\beta$ on $(y_{\ell},y_r)$
such that
 \[
  \dspt{\Psi(y)=\frac{1}{2m}\,\beta(y)^{2m}.}
 \]
We take this diffeomorphism into account and we can write the
involution associated to $\Psi$ as
\begin{equation*}
  \sigma_2(y)=\beta^{-1}\bigl(-\beta(y)\bigr).
\end{equation*}
In what follows, for each $h\in (0,h_0)$, we denote the projection
of the oval $\gamma_h$ on the $x$-axis by $(x_h^-,x_h^+).$
Therefore, $x_{\ell}<x_h^-<0<x_h^+<x_r$ and $\Phi(x_h^{\pm})=h.$
 \begin{figure}[t]
  \begin{center}
   \epsfig{file=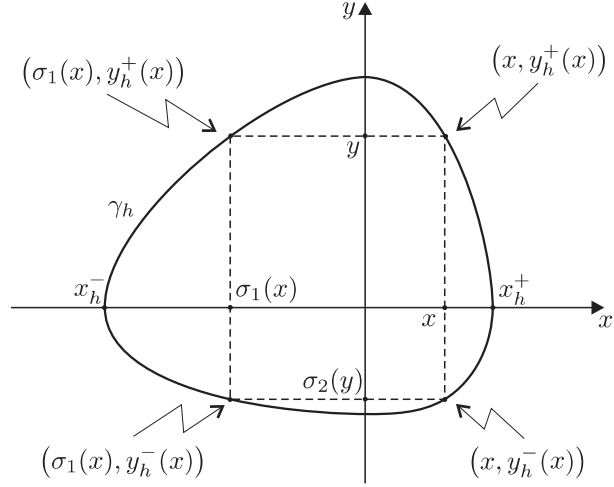}
  \end{center}
  \caption{Notation related to the oval $\gamma_h$.}
  \label{dib1}
 \end{figure}
Moreover (see \figc{dib1}), if $(x,y)\in\gamma_h,$ then
 \[
  y= y_h^+(x)\,\mbox{ for $y>0$ and }y= y_h^-(x)\,\mbox{ for $y<0,$}
 \]
where
 \[
 y_h^{\pm}(x)\!:=\beta^{-1}\left(\pm\sqrt[2m]{2m\bigl(h-\Phi(x)\bigr)}\right).
 \]
We note that $y_{h}^{\pm}(x)=y_{h}^{\pm}\bigl(\sigma_1(x)\bigr),$
where we recall that $\sigma_1$ is the involution associated to
$\Phi.$ We begin by the proof of the following result.

\begin{lem}\label{lem2}
Let $f$ and $g$ be analytic functions on $(x_{\ell},x_r)$ and
$(y_{\ell},y_r)$, respectively, and let us consider
 \[
 I(h)=\int_{\gamma_h}f(x)g(y)dx.
 \]
We set $\ell(x)\!:= f(x)-f\bigl(\sigma_1(x)\bigr)\sigma'_1(x)$ and
$\xi_k\!:=\mathscr B_{\sigma_2}(g_k),$ where $g_k$ is recursively
defined by means of $g_{k+1}=\frac{g'_{k}}{\Psi'}$ with $g_0= g.$
Then, if $\xi_0(y)=\op\bigl(y^{2m(n-2)}\bigr),$
 \[
  I^{(k)}(h)=\int_0^{x_h^+}\ell(x)\,\xi_k\bigl(y_h^+(x)\bigr)dx\
  \mbox{ for $k=0,1,\ldots,n-1.$}
 \]
\end{lem}

\begin{prova} We prove the result by induction on $k.$
We take the parameterization of the oval $\gamma_h$ given by the
mappings $x\longmapsto\bigl(x,y_h^{\pm}(x)\bigr)$, with the
clockwise orientation, and we use
$y_h^-(x)=\sigma_2\bigl(y_h^+(x)\bigr),$ to get that
 \begin{align*}
  I(h)=&\int_{x_h^+}^{x_h^-}f(x)g\bigl(y_h^-(x)\bigr)dx+\int_{x_h^-}^{x_h^+}f(x)g\bigl(y_h^+(x)\bigr)dx
  =\int_{x_h^-}^{x_h^+}f(x)\left.\bigl(g(y)-g(\sigma_2(y))\bigr)\right|_{y=y_h^+(x)}dx \\[5pt]
  =&\int_{x_h^-}^{0}f(x)\left.\bigl(g(y)-g(\sigma_2(y))\bigr)\right|_{y=y_h^+(x)}dx +
  \int_{0}^{x_h^+}f(x)\left.\bigl(g(y)-g(\sigma_2(y))\bigr)\right|_{y=y_h^+(x)}dx\\[5pt]
  =&\int_{x_h^+}^{0}f\bigl(\sigma_1(u)\bigr)\sigma_1'(u)
  \left.\bigl(g(y)-g(\sigma_2(y))\bigr)\right|_{y=y_h^+(\sigma_1(u))}du+
  \int_{0}^{x_h^+}f(x)\left.\bigl(g(y)-g(\sigma_2(y))\bigr)\right|_{y=y_h^+(x)}dx,
 \end{align*}
where in the last equality we performed the change of variable
$x=\sigma_1(u).$ Thus, since $y_h^+(\sigma_1(u))=y_h^+(u),$ the
above expression yields to
 \[
  I(h)=\int_{0}^{x_h^+}\bigl(f(x)-f(\sigma_1(x))\sigma_1'(x)\bigr)  \left. \bigl(g(y)
  -g(\sigma_2(y))\bigr)\right|_{y=y_h^+(x)}dx=\int_{0}^{x_h^+}\ell(x)
  \mathscr B_{\sigma_2}\bigl(g\bigr)\bigl(y_h^+(x)\bigr)dx.
 \]
This expression proves the result for $k=0.$  We assume now that
the result holds true for $k<n-1.$ On account of the hypothesis
about the order of $\xi_0$ at $y=0,$ an easy computation shows
that $\xi_k(y)=\mathscr
B_{\sigma_2}\bigl(g_k\bigr)(y)=\op\bigl(y^{2m(n-2-k)}\bigr).$ The
fact that ${2m(n-2-k)}\geqslant 0$ enables us to differentiate the
expression of $I^{(k)}(h)$ and we obtain
 \begin{align*}
  I^{(k+1)}(h)&=\frac{d}{dh}\int_0^{x_h^+}\ell(x)\xi_k\bigl(y_h^+(x)\bigr)dx\\[5pt]
  &=\ell\bigl(x_h^+\bigr)\xi_k(0)\frac{dx_h^+(x)}{dh}
  +\int_0^{x_h^+}\ell(x)\xi_k'\bigl(y_h^+(x)\bigr)\frac{dy_h^+(x)}{dh}dx
  =\int_0^{x_h^+}\ell(x)\left.\frac{\xi_k'(y)}{\Psi'(y)}\right|_{y=y_h^+(x)}dx
 \end{align*}
(Let us note that in the second equality we use that $y^+_h(x)=0$
at $x=x_h^+$ because $\Phi(x_h^+)=h$ and
$\Psi\bigl(y_h^+(x)\bigr)=h$ for all $h.)$ Finally, since
 \[
  \dspt{\displaystyle \xi_k'(y)=g_k'(y)-g_k'\bigl(\sigma_2(y)\bigl)\sigma'_2(y)
  =g_k'(y)-g_k'\bigl(\sigma_2(y)\bigl)\frac{\Psi'(y)}{\Psi'\bigl(\sigma_2(y)\bigr)}
  =\Psi'(y)\mathscr B_{\sigma_2}\bigl(\frac{g_k'}{\Psi'}\bigr)(y)
  =\Psi'(y)\xi_{k+1}(y),}
 \]
the result for $k+1$ follows and the proof is completed.
\end{prova}

\begin{obs}\label{impropia}
It is worth making some comments on the expression of the $(n-1)$
derivative of $I(h)$ given by \lemc{lem2}. The condition $\mathscr
B_{\sigma_2}(g_0)(y)=\xi_0(y)=\op\bigl(y^{2m(n-2)}\bigr)$
guarantees that the integral
 \[
 I^{(n-1)}(h)=\int_0^{x_h^+}\ell(x)\,\xi_{n-1}\bigl(y_h^+(x)\bigr)dx,
 \]
despite it may be improper, is convergent. Indeed, by this
condition, the Taylor series of $\xi_0$ at $y=0$ begins at least
with order $2m(n-2)+1,$ i.e. $\xi_0(y)=\Delta
y^{2m(n-2)+1}+\ldots$ with $\Delta\neq 0.$ To construct
$g_{k+1}(y)$, we derive $g_k(y)$ and divide it by $\Psi'(y),$
which vanishes at $y=0$ with multiplicity $2m-1.$ Hence, it turns
out that $\xi_{n-1}=\mathscr B_{\sigma_2}(g_{n-1})$ is not
analytic at $y=0$ but meromorphic. However, due to the mentioned
condition, the pole has at most order~$2m-1.$ We note that
$y_h^+(x)=0$ at $x=x_h^+$ because $\Phi(x_h^+)=h.$ More precisely,
we take $\Phi'(x_h^+)\neq 0$ also into account and it is easy to
show that
 \[
  \lim_{x\longrightarrow
  x_h^+}\dspt{\frac{y_h^+(x)}{\sqrt[2m]{x-x_h^+}}}\neq 0.
 \]
Accordingly, although $\xi_{n-1}\bigl(y_h^+(x)\bigr)$ may tend to
infinity as $x\longrightarrow x_h^+,$ the derivative $I^{(n-1)}(h)$
is given by a convergent integral.
\end{obs}

Let us consider now
 \begin{equation*}
  I_k(h)=\int_{\gamma_h}f_k(x)g(y)dx,\,\mbox{ for $k=0,1,\ldots,n-1,$}
 \end{equation*}
where $g$ is an analytic function on $(y_{\ell},y_r)$ and each
$f_k$ is an analytic function on $(x_{\ell},x_r).$ The next result
provides an expression of the Wronskian of
$(I_0,I_1,\ldots,I_{k-1}).$ In its statement, $\xi_i$ is defined
as in \lemc{lem2}, i.e. we set $g_{i+1}=\frac{g'_{i}}{\Psi'}$ with
$g_0=g,$ and $\xi_i\!:=\mathscr B_{\sigma_2}(g_i).$ Moreover
 \[
  \Delta_k(h)\!:=\bigl\{\mathbf{x_k}\in\R^k:0<x_0<x_1<\ldots<x_{k-1}<x_h^+\bigr\}.
 \]

\begin{prop}\label{prop3}
Let us assume that  $\mathscr
B_{\sigma_2}\bigl(g\bigr)(y)=\op\bigl(y^{2m(n-2)}\bigr).$ Then,
for each $k=1,2,\ldots,n,$ the Wronskian of
$(I_0,I_1,\ldots,I_{k-1})$ at $h\in (0,h_0)$ is given by
 \[
   W\bigl[\mathbf{I_k}\bigr](h)=\int\cdots\int_{\Delta_k(h)}D\bigl[\mathbf{\ell_k}\bigr](\mathbf{x_k})
   D\bigl[\mathbf{\xi_k}\bigr](\mathbf{y_k})\,dx_0\,dx_1\cdots dx_{k-1},
 \]
where $y_i=y_h^+(x_i)$ and
$\ell_i(x)=f_i(x)-f_i\bigl(\sigma_1(x)\bigr)\sigma_1'(x).$
\end{prop}

\begin{prova}
Fix $k\in\{1,2,\ldots,n\}$ and let $S_k$ be the symmetric group of
$k$ elements. We take the definition of determinant into account and
we apply \lemc{lem2} to show that
 \begin{align*}
 W\bigl[\mathbf{I_k}\bigr](h)=&\det\left(I_j^{(i)}(h)\right)_{0\leqslant i,j\leqslant k-1}
 =\sum_{\tau\in S_k}\mbox{sgn}(\tau)\prod_{i=0}^{k-1}I_{\tau(i)}^{(i)}(h) \\[5pt]
 =&\sum_{\tau\in S_k}\mbox{sgn}(\tau)\prod_{i=0}^{k-1}\int_0^{x_h}\ell_{\tau(i)}(x)\,\xi_i\bigl(y_h^+(x)\bigr)dx\\[5pt]
 =&\sum_{\tau\in S_k}\mbox{sgn}(\tau)\prod_{i=0}^{k-1}\int_0^{x_h}\ell_{\tau(i)}(x_i)\,\xi_i\bigl(y_h^+(x_i)\bigr)dx_i
 \\[5pt]
 =&\int\cdots\int_{[0,x_h^+]^k}\left[\sum_{\tau\in S_k}\mbox{sgn}(\tau)
   \prod_{i=0}^{k-1}\ell_{\tau(i)}(x_i)\right]\prod_{i=0}^{k-1}\xi_i(y_i)\,dx_0\,dx_1\cdots dx_{k-1}\\[5pt]
 =&\int\cdots\int_{[0,x_h^+]^k}D\bigl[\mathbf{\ell_k}\bigr](\mathbf{x_k})
   \prod_{i=0}^{k-1}\xi_i(y_i)\,dx_0\,dx_1\cdots dx_{k-1}.
 \end{align*}
At this point, for each permutation $\tau\in S_k$ we define
\map{\psi_{\tau}}{\R^k}{\R^k} as
 \[
  \psi_{\tau}(x_0,x_1,\ldots,x_{k-1})=(x_{\tau(0)},x_{\tau(1)},\cdots,x_{\tau(k-1)}),
 \]
which is clearly an invertible mapping. We note that
 \[
  [0,x_h^+]^k\setminus\mathcal R=\bigcup_{\tau\in
  S_k}\psi_{\tau}\bigl(\Delta_k(h)\bigr),
 \]
where $\mathcal R$ is a subset of $\R^k$ with Lebesgue measure equal
to zero. Accordingly
 \begin{align*}
  W\bigl[\mathbf{I_k}\bigr](h)=&\int\cdots\int_{[0,x_h^+]^k}D\bigl[\mathbf{\ell_k}\bigr](\mathbf{x_k})
    \prod_{i=0}^{k-1}\xi_i(y_i)\,dx_0\,dx_1\cdots dx_{k-1} \\[5pt]
    =&\sum_{\tau\in S_k}\int\cdots\int_{\psi_{\tau}\left(\Delta_k(h)\right)}
    D\bigl[\mathbf{\ell_k}\bigr](\mathbf{x_k})\prod_{i=0}^{k-1}\xi_i(y_i)\,dx_0\,dx_1\cdots dx_{k-1}.
 \end{align*}
Next, in each integral of the above summation we perform the
coordinate transformation $\mathbf{x_k}=\psi_{\tau}(\mathbf{u_k})$
(i.e., $x_i=u_{\tau(i)}$ for $i=0,1,\ldots,k-1$), so that
 \[
  W\bigl[\mathbf{I_k}\bigr](h)
    =\sum_{\tau\in S_k}\int\cdots\int_{\Delta_k(h)}
    D\bigl[\mathbf{\ell_k}\bigr]\bigl(\psi_{\tau}(\mathbf{u_k})\bigr)\prod_{i=0}^{k-1}
     \xi_i\bigl(v_{\tau(i)}\bigr)\,du_0\,du_1\cdots du_{k-1},
 \]
where $v_i=y_h^+(u_i).$ (Here we use that the absolute value of
the determinant of the Jacobian of~$\psi_{\tau}$ is identically
one.) Finally, we remark that
$D\bigl[\mathbf{\ell_k}\bigr]\bigl(\psi_{\tau}(\mathbf{u_k})\bigr)
\, = \,
\mbox{sgn}(\tau)D\bigl[\mathbf{\ell_k}\bigr](\mathbf{u_k})$ and we
take the properties of the determinant into account to prove that
 \begin{align*}
  W\bigl[\mathbf{I_k}\bigr](h)
    &=\sum_{\tau\in S_k}\int\cdots\int_{\Delta_k(h)}
    \mbox{sgn}(\tau)D\bigl[\mathbf{\ell_k}\bigr](\mathbf{u_k})\prod_{i=0}^{k-1}\xi_i\bigl(v_{\tau(i)}\bigr)\,du_0\,du_1\cdots
    du_{k-1} \\[5pt]
    &=\int\cdots\int_{\Delta_k(h)}D\bigl[\mathbf{\ell_k}\bigr](\mathbf{u_k})
    \left(\sum_{\tau\in S_k}\mbox{sgn}(\tau)
    \prod_{i=0}^{k-1}\xi_i\bigl(v_{\tau(i)}\bigr)\right)du_0\,du_1\cdots du_{k-1}\\[5pt]
    &=\int\cdots\int_{\Delta_k(h)}D\bigl[\mathbf{\ell_k}\bigr](\mathbf{u_k})
    D[\mathbf{\xi_k}](\mathbf{v_k})\,du_0\,du_1\cdots du_{k-1},
 \end{align*}
and this last identity proves the result.
\end{prova}

\begin{prooftext}{Proof of \teoc{desacoplat}.}
We claim that the assumptions $(a)$ and $(b)$ imply that the
Wronskians $W\bigl[\mathbf{I_k}\bigr](h)$ for $k=1,2,\ldots,n$ are
different from zero at any $h\in (0,h_0).$ On account of $(b)$ in
\lemc{lem1}, this fact will prove that $(I_0,I_1,\ldots,I_{n-1})$
is an ECT-system on $(0,h_0).$

$ $From \propc{prop3},
 \[
  W\bigl[\mathbf{I_k}\bigr](h)=\int\cdots\int_{\Delta_k(h)}D\bigl[\mathbf{\ell_k}\bigr](\mathbf{x_k})
    D\bigl[\mathbf{\xi_k}\bigr](\mathbf{y_k})\,dx_0\,dx_1\cdots dx_{k-1},
 \]
where recall that
$y_i=y_h^+(x_i)=\beta^{-1}\left(\sqrt[2m]{2m(h-\Phi(x_i))}\right)$.
On the other hand, $x\longmapsto
\beta^{-1}\left(\sqrt[2m]{2m(h-\Phi(x_i))}\right)$ is decreasing on
$(0,x_r)$ and, therefore, in the above integral we have that
 \[
  0<x_0<x_1<\ldots<x_{k-1}< x_h^+\mbox{ and }0<y_{k-1}<y_{k-2}<\ldots<y_0<y_h^+.
 \]
We note at this point that $\ell_i(x)=\Phi'(x)\mathscr
B_{\sigma_1}\!\bigl(\frac{f_i}{\Phi'}\bigr)(x)$ because
\[
  \dspt{\ell_i(x)=f_i(x)-f_i\bigl(\sigma_1(x)\bigr)\sigma_1'(x)=f_i(x)-f_i\bigl(\sigma_1(x)\bigr)\frac{\Phi'(x)}
  {\Phi'\bigl(\sigma_1'(x)\bigr)}=\Phi'(x)\left(\bigl(\frac{f_i}{\Phi'}\bigr)(x)
  -\bigl(\frac{f_i}{\Phi'}\bigr)\bigl(\sigma_1(x)\bigr)\right)}.
\]
Since $\Phi'(x) \neq 0$ for any $x \in \ii$ and, by assumption,
$\Bigl(\mathscr
B_{\sigma_1}\!\bigl(\frac{f_0}{\Phi'}\bigr),\mathscr
B_{\sigma_1}\!\bigl(\frac{f_1}{\Phi'}\bigr),\ldots,\mathscr
B_{\sigma_1}\!\bigl(\frac{f_{n-1}}{\Phi'}\bigr)\Bigr)$ is a
CT-system on $(0,x_r),$ so it is
$(\ell_0,\ell_1,\ldots,\ell_{n-1}).$ The second assumption ensures
that $(\xi_0,\xi_1,\ldots,\xi_{n-1})$ is a CT-system on $(0,y_r)$
because, by definition, $\xi_i=\mathscr B_{\sigma_2}(g_i).$
Therefore, we apply statement $(a)$ in \lemc{lem1} and it turns
out that
 \[
  D\bigl[\mathbf{\ell_k}\bigr](\mathbf{x_k})D\bigl[\mathbf{\xi_k}\bigr](\mathbf{y_k})\neq 0
  \mbox{ for all $\mathbf{x_k}\in\Delta_k(h).$}
 \]
Since $\Delta_k(h)$ is connected, we have shown that
$W\bigl[\mathbf{I_k}\bigr](h)\neq 0$ and the result follows.
\end{prooftext}

\begin{prooftext}{Proof of \teoc{quadratic}.}
This result is in fact a corollary of \teoc{desacoplat}. We note
that $B(x)>0$ for $x\in (x_\ell, x_r)$. Thus the coordinate
transformation
$(u,v)=\chi(x,y)\!:=\bigl(x,\sqrt[2m]{2mB(x)}\,y\bigr)$ is well
defined and verifies
$e_h\!:=\chi^{-1}(\gamma_h)\subset\bigl\{A(u)+\frac{1}{2m}v^{2m}=h\bigr\}.$
Accordingly
 \[
  I_i(h)=\int_{\gamma_h}f_i(x)y^{2s-1}dx=(2m)^{\frac{1-2s}{2m}}\int_{e_h}
    \dspt{\left(\frac{f_i}{B^{\frac{2s-1}{2m}}}\right)}\!(u)\,v^{2s-1}du.
 \]
Following the obvious notation, we can apply \teoc{desacoplat}
with
 \[
  \dspt{\widehat f_i=\frac{f_i}{B^{\frac{2s-1}{2m}}},\quad\widehat g(v)=v^{2s-1},\quad\Phi=A,\quad
 \Psi(v)=\frac{1}{2m}\,v^{2m},\quad\sigma_1=\sigma\;\mbox{ and }\;\sigma_2=-Id.}
 \]
Clearly the hypothesis $(a)$ in \teoc{desacoplat} is guaranteed by
the assumption on $\ell_i=\mathscr
B_{\sigma_1}\!\left(\frac{\widehat f_i}{\Phi'}\right).$ Let us turn
now to the hypothesis $(b).$ We take $\sigma_2=-Id$ and
$\Psi'(v)=v^{2m-1}$ into account and one can easily show that
$\widehat g_i(v)=c_iv^{2(s-im)-1}$ for some positive constant $c_i,$
so that $\mathscr B_{\sigma_2}(\widehat g_i)(v)=2c_iv^{2(s-im)-1}.$
Hence, $\Bigl(\mathscr B_{\sigma_2}\!(\widehat g_0),\mathscr
B_{\sigma_2}\!(\widehat g_1),\ldots,\mathscr B_{\sigma_2}\!(\widehat
g_{n-1})\Bigr)$ is clearly a CT-system on $(0,+\infty).$ Since the
condition $s>m(n-2)$ implies that $\mathscr B_{\sigma_2}(\widehat
g)(v)=2v^{2s-1}=\op(v^{2m(n-2)}),$ the hypothesis $(b)$ in
\teoc{desacoplat} is satisfied as well. Therefore, we apply
\teoc{desacoplat} and we can assert that $(I_0,I_1,\ldots,I_{n-1})$
is an ECT-system on $(0,h_0)$ as desired.
\end{prooftext}

\section{Applications}\label{Sec4}

The following lemma establishes a formula to write the integrand of
an Abelian integral so as to be suitable to apply our results.

\begin{lem}\label{puja}
Let $\gamma_h$ be an oval inside the level curve
$\{A(x)+B(x)y^2=h\}$ and we consider a function $F$ such that
$F/A'$ is analytic at $x=0.$ Then, for any $k\in\N,$
 \[
  \int_{\gamma_h}F(x)y^{k-2}dx=\int_{\gamma_h}G(x)y^kdx
 \]
where
$G(x)=\frac{2}{k}\bigl(\frac{BF}{A'}\bigr)'\!(x)-\bigl(\frac{B'F}{A'}\bigr)(x).$
\end{lem}

\begin{prova}
If $(x,y)\in\gamma_{h}\subset\{A(x)+B(x)y^2=h\}$ then
$\frac{dy}{dx}=-\frac{A'(x)+B'(x)y^2}{2B(x)y},$ and accordingly
 \begin{align*}
  d\bigl(g(x)y^k\bigr)=&g'(x)y^kdx+kg(x)y^{k-1}dy \\[3pt]
  =&
  \dspt{\left(g'(x)-\frac{k}{2}\bigl(\frac{A'g}{B}\bigr)(x)\right)y^kdx
  -\frac{k}{2}\bigl(\frac{A'g}{B}\bigr)(x)\,y^{k-2}dx.}
 \end{align*}
We take $F(x)=\frac{k}{2}\bigl(\frac{A'g}{B}\bigr)(x)$ in the
above equality, we use that
$\int_{\gamma_h}d\bigl(g(x)y^k\bigr)=0$ and the result follows.
\end{prova}

$ $From now on we shall often compute the resultant between two
polynomials and we shall apply Sturm's Theorem to study the number
of roots of a polynomial in an interval. The interested reader is
referred to the Appendix for details.

\begin{ex}\label{ex1}
Iliev and Perko study in~\cite{Iliev1} symmetric Hamiltonian
systems perturbed asymmetrically. More concretely, systems of the
form
 \begin{equation*}
  \left\{
   \begin{array}{l}
    \dot x=y, \\[4pt]
    \dot y=\pm(x\pm x^3)+\lambda_1 y+\lambda_2 x^2+\lambda_3 xy+\lambda_4 x^2y,
   \end{array}
  \right.
 \end{equation*}
where $\lambda_j(\varepsilon)=O(\varepsilon),$ and they prove that
at most two limit cycles bifurcate for small $\varepsilon\neq 0$
from any period annulus of the unperturbed system. There are three
different cases to consider depending on the phase portrait of the
unperturbed system: the \emph{global center}, the \emph{truncated
pendulum} and the \emph{Duffing oscillator.} This latter case gives
rise to two different types of period annuli (see \figc{dib2}).
 \begin{figure}[t]
  \begin{center}
   \epsfig{file=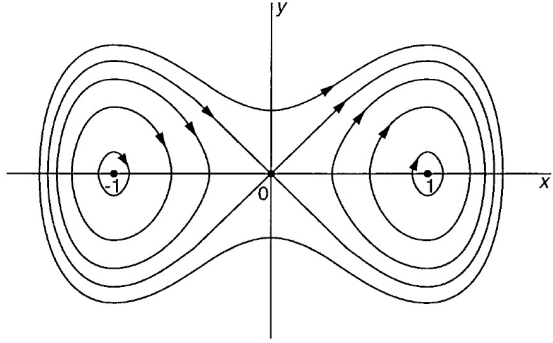,scale=0.8}
  \end{center}
  \caption{The period annuli in the Duffing oscillator.}
  \label{dib2}
 \end{figure}
In this example we study the so-called \emph{interior Duffing
oscillator.} Theorem 1.3 in~\cite{Iliev1} shows that at most two
limit cycles bifurcate from either one of the interior period
annuli.

If we perform a translation to bring the center on the right
half-plane to the origin, the Hamiltonian function of the
unperturbed system becomes
 \[
  H(x,y)=A(x)+B(x)y^2\,\mbox{ with $A(x)=x^2+\frac{1}{4}x^4+x^3$ and $B(x)=\frac{1}{2}.$}
 \]
The projection of the period annulus of this center is
$\bigl(-1,\sqrt{2}-1\bigr)$ and $h_0=A(-1)=1/4.$

$ $From Theorem~2.1 in~\cite{Iliev1}, it follows that the first
non-identically zero Melnikov function is a linear combination of
$\widetilde I_i(h)=\int_{\gamma_h}x^iydx$ for $i=0,1,2.$ Thus,
Theorem~1.3 in~\cite{Iliev1} will follow if we prove that
$\bigl\{\widetilde I_0,\widetilde I_1,\widetilde I_2\bigr\}$ is an
ECT-system. Additionally, this fact implies that there are values
of the parameters for which exactly $0$, $1$ or $2$ limit cycles
bifurcate from the period annulus. To this end we will apply
\teoc{quadratic}, but we note that in this case $m=1,$ $n=3$ and
$s=1,$ so that the hypothesis $s>m(n-2)$ is not satisfied. This is
easy to overcome because
 \[
  \widetilde I_0(h)=\int_{\gamma_h}ydx=\frac{1}{h}\int_{\gamma_h}\bigl(A(x)+B(x)y^2\bigr)ydx
  =\frac{1}{h}\int_{\gamma_h}A(x)ydx+\frac{1}{h}\int_{\gamma_h}B(x)y^3dx,
 \]
and then, we apply \lemc{puja} with $k=3$ and $F=A$ to the first
integral above, to get
 \[
  \widetilde I_0(h)=\frac{1}{h}\int_{\gamma_h}\frac{x^2+2x+2}{12(x+1)^2}\,y^3dx
    +\frac{1}{h}\int_{\gamma_h}\frac{1}{2}\,y^3dx
    =\frac{1}{h}\int_{\gamma_h}f_0(x)y^3dx\,
    \mbox{ with }f_0(x)\!:=\frac{7x^2+14x+8}{12(x+1)^2}.
 \]
(It is not possible to apply \lemc{puja} directly to $\widetilde
I_0$ because then we must take $F\equiv 1,$ and in this case $F/A'$
is not analytic at $x=0.)$ Exactly in the same way we obtain
 \[
 \begin{array}{ll}
  \dsp \widetilde I_1(h)=\frac{1}{h}\int_{\gamma_h}f_1(x)y^3dx&\,\mbox{ with
  }\dsp f_1(x)\!:=\frac{x(8x^2+17x+10)}{12(x+1)^2}, \\[15pt]
  \dsp \widetilde I_2(h)=\frac{1}{h}\int_{\gamma_h}f_2(x)y^3dx&\,\mbox{ with
  }\dsp f_2(x)\!:=\frac{x^2(9x^2+20x+12)}{12(x+1)^2}.
 \end{array}
 \]
We set $I_i(h)=\int_{\gamma_h}f_i(x)y^3dx$ and it is clear that
$\bigl\{\widetilde I_0,\widetilde I_1,\widetilde I_2\bigr\}$ is an
ECT-system on $(0,h_0)$ if and only if so it is $\{I_0,I_1,I_2\}.$
We can now apply \teoc{quadratic} because $s=2$ and the condition
$s>m(n-2)$ holds. Thus, setting
 \[
  \dspt{\ell_i(x)=\left(\frac{f_i}{A'}\right)\!(x)
   -\left(\frac{f_i}{A'}\right)\!\bigl(\sigma(x)\bigr)},
 \]
we have to check that $\{\ell_0,\ell_1,\ell_2\}$ is a CT-system on
$\bigl(0,\sqrt{2}-1\bigr).$ Here $\sigma$ is the involution
associated to~$A$ and we used that~$B$ is constant. (In this
example we can compute the involution explicitly but we do not use
it because we want to show that it is not necessary to apply our
result.) As a matter of fact we will show that
$\{\ell_0,\ell_1,\ell_2\}$ is an ECT-system because a continuous
Wronskian is easy to study. In order to compute the three
Wronskians, we write $\ell_i(x)=L_i\bigl(x,\sigma(x)\bigr)$ with
$L_i(x,z)=\bigl(\frac{f_i}{A'}\bigr)(x)
-\bigl(\frac{f_i}{A'}\bigr)(z).$ Moreover, due to
 \[
  \dspt{A(x)-A(z)=\frac{1}{4}(x-z)(x+2+z)(x^2+2x+2z+z^2)},
 \]
it turns out that $z=\sigma(x)$ is defined by means of
$q(x,z)\!:=x^2+2x+2z+z^2= 0.$ Accordingly, since
$\sigma'(x)=-\frac{x+1}{z+1},$ we have that
$W[\mathbf{\ell_i}](x)=\omega_i\bigl(x,\sigma(x)\bigr)$ with
$\omega_i(x,z)$ being a \emph{rational} function for $i=1,2,3.$ The
resultant with respect to $z$ between $q(x,z)$ and the numerator of
$\omega_3(x,z)$ is $r_3(x)=64x^{16}(x+2)^{16}p_3(x)$ with
 \begin{align*}
  p_3(x)=&\,441\,{x}^{20}+8820\,{x}^{19}+79380\,{x}^{18}+423360\,{x}^{17}+1481685
      \,{x}^{16}+3555024\,{x}^{15}+5918640\,{x}^{14}\\&+6740160\,{x}^{13}+
      4976155\,{x}^{12}+1881540\,{x}^{11}-892716\,{x}^{10}-3303200\,{x}^{9}-
      4779945\,{x}^{8}\\&-3240840\,{x}^{7}+601960\,{x}^{6}+2523360\,{x}^{5}+
      1158080\,{x}^{4}-414400\,{x}^{3}-414400\,{x}^{2}+44800,
 \end{align*}
and by applying Sturm's Theorem we can assert that $p_3(x)\neq 0$
for all $x\in\bigl(0,\sqrt{2}-1\bigr).$ Thus, $\omega_3(x,z)=0$ and
$q(x,z)=0$ have no common roots, and this fact implies that
$W[\mathbf{\ell_3}](x)\neq 0$ for all
$x\in\bigl(0,\sqrt{2}-1\bigr).$ The resultant with respect to $z$
between $q(x,z)$ and the numerator of $\omega_2(x,z)$ is
$r_2(x)=32x^7(x+2)^7p_2(x)$ with
 \begin{align*}
  p_2(x)=&\,49\,{x}^{12}+588\,{x}^{11}+2940\,{x}^{10}+7840\,{x}^{9}+11650\,{x}^{8}
         +8528\,{x}^{7}\\&+496\,{x}^{6}-3520\,{x}^{5}-1915\,{x}^{4}-620\,{x}^{3}-
         620\,{x}^{2}+360,
 \end{align*}
and using Sturm's Theorem it follows that $p_2$ does not vanish on
$(0,\sqrt{2}-1).$ Exactly as before, this fact shows that
$W[\mathbf{\ell_2}](x)\neq 0$ for all
$x\in\bigl(0,\sqrt{2}-1\bigr).$ Finally, the resultant with
respect to $z$ between $q(x,z)$ and the numerator of
$\omega_1(x,z)$ is
 \[
  r_1(x)=2\,{x}^{3} \left( x+2 \right) ^{3} \left( 49\,{x}^{8}+392\,{x}^{7}+
1176\,{x}^{6}+1568\,{x}^{5}+659\,{x}^{4}-500\,{x}^{3}-500\,{x}^{2}+80\right)
 \]
and, thanks to Sturm's Theorem again, we can assert that it does
not vanish on $(0,\sqrt{2}-1).$ This proves that
$W[\mathbf{\ell_1}](x)=\ell_0(x)\neq 0$  for all $x\in
(0,\sqrt{2}-1).$ Consequently $\{\ell_0,\ell_1,\ell_2\}$ is an
ECT-system on $(0,\sqrt{2}-1)$ and by applying \teoc{quadratic},
$\{I_0,I_1,I_2\}$ is an ECT-system on $(0,1/4).$ Therefore, the
first Melnikov function has at most two zeros counting
multiplicities.
\end{ex}

\begin{ex}\label{ex2}
Zhao, Liang and Lu study in \cite{Zhao} the system of planar
differential equations
 \begin{equation*}
  \left\{
   \begin{array}{l}
    \dot x=2xy+\varepsilon\Bigl(\dsp\sum_{i+j\leqslant 2}a_{ij}(\varepsilon)x^iy^j\Bigr), \\[12pt]
    \dot y=6x-6x^2-y^2+\varepsilon\Bigl(\dsp\sum_{i+j\leqslant
    2}b_{ij}(\varepsilon)x^iy^j\Bigr).
   \end{array}
  \right.
 \end{equation*}
The unperturbed system (i.e., with $\varepsilon =0$) has a center
at $(1,0)$ whose period annulus is bounded by a cuspidal loop and
they prove (see Theorem~1.2 in~\cite{Zhao}) that the maximum
number of limit cycles emerging from its period annulus for
$\varepsilon\approx 0$ is two.

Our goal is to reobtain this result by applying \teoc{quadratic}. To
this end, we bring the center to the origin by means of a
translation, so that the unperturbed system is Hamiltonian with
 \[
  H(x,y)=A(x)+B(x)y^2,\,\mbox{ where $A(x)=x^2(3+2x)$ and $B(x)=x+1.$}
 \]
The projection of the period annulus is now $(-1,1/2)$ and the
energy level of the polycycle in its outer boundary is
$h_0=A(-1)=1.$ By Theorem~3 in~\cite{Iliev2}, the upper bound for
the number of limit cycles is equal to the maximum number of zeros
for $h\in (0,1),$ counted with multiplicities, of any non-trivial
linear combination of
 \[
  \widetilde I_i(h)=\int_{\gamma_h}(x+1)^{i-1}ydx\,\mbox{ for $i=0,1,2.$}
 \]
Accordingly, the result in~\cite{Zhao} will follow once we show
that $\bigl\{\widetilde I_0,\widetilde I_1,\widetilde I_2\bigr\}$
is an ECT-system on $(0,1).$ By applying \lemc{puja}, the same
straightforward manipulation as before shows that $\widetilde
I_i(h)=\frac{1}{18h}I_i(h)$ where
 \[
 I_i(h)=\int_{\gamma_h}f_i(x)y^3dx
 \]
with
 \[
  f_0(x)=\frac{16x^2+35x+24}{(x+1)^2}\,,\ \mbox{ }f_1(x)=\frac{20x^2+41x+24}{x+1}\, \ \mbox{ and
  }\   f_2(x)=24x^2+47x+24.
 \]
It is clear that $\bigl\{\widetilde I_0,\widetilde I_1,\widetilde
I_2\bigr\}$ is an ECT-system on the interval $(0,1)$ if, and only
if, so it is $\{I_0,I_1,I_2\}.$ On account of \teoc{quadratic}, this
will follow once we check that $\{\ell_0,\ell_1,\ell_2\}$ is an
ECT-system on $(0,1/2),$ where $\ell_i=\mathscr
B_{\sigma}\!\left(\frac{f_i}{A'B^{3/2}}\right).$ Note that
$A(x)-A(z)=(x-z)(2x^2+2zx+3x+2z^2+3z),$ so that $z=\sigma(x)$ is
implicitly defined by means of $q(x,z)\!:=2x^2+2zx+3x+2z^2+3z=0.$
Thus
 \[
  \sigma'(x)=\frac{dz}{dx}=-\frac{4x+2z+3}{4z+2x+3}.
 \]
Taking this into account, some computations show that, for
$i=1,2,3,$ $W[\,\mathbf{\ell_i}](x)=\omega_i\bigl(x,\sigma(x)\bigr)$
with $\omega_i(x,z)$ being a \emph{rational} function of
$u=\sqrt{x+1}$ and $v=\sqrt{z+1},$ say $R_i(u,v).$ Note that
$x\longmapsto\sqrt{x+1}$ maps $(0,1/2)$ to $(1,\sqrt{3/2}).$ The
resultant with respect to $v$ between the numerator of $R_i(u,v)$
and $q(u^2-1,v^2-1)$ is a polynomial $r_i(u)$ that, by applying
Sturm's Theorem, has no roots on $(1,\sqrt{3/2}).$ (For the sake of
shortness we do not give here the expression of these polynomials.)
Hence, it is proved that $W[\,\mathbf{\ell_i}]$ does not vanish on
$(0,1/2)$ for $i=1,2,3.$ By \teoc{quadratic}, this reasoning proves
the mentioned result of Zhao, Liang and Lu.
\end{ex}

\begin{ex}\label{ex3}
Peng studies in \cite{Peng} the system of planar differential
equations
 \begin{equation*}
  \left\{
   \begin{array}{l}
    \dot x=-y-3x^2-y^2+\varepsilon(\mu_1x+\mu_2xy), \\[5pt]
    \dot y=x(1-2y)+\varepsilon\mu_3x^2.
   \end{array}
  \right.
 \end{equation*}
The unperturbed system (i.e. when $\varepsilon =0$) has a center
at the origin and the author proves (see Theorem~A in~\cite{Peng})
that two is the maximal number of limit cycles which bifurcate
from its period annulus for $\varepsilon\approx 0$ and that there
are perturbations with exactly $0$, $1$ or $2$ limit cycles. To
this end, he first shows that by means of the projective
coordinate transformation $(x,y) \mapsto
(\frac{y}{x+2},\frac{x}{2(x+2)})$ and a non-constant rescaling of
time the above system reads for
 \begin{equation*}
  \left\{
   \begin{array}{l}
    \dot x=2(x+2)y + \varepsilon \mu_3 (x+2)y^2,\\[5pt]
    \dot y=-x-\frac{3}{4} x^2-y^2+\varepsilon \bigl( \mu_1 (x+2) + \frac{\mu_2}{2} x + \mu_3 y^2 \bigr).
   \end{array}
  \right.
 \end{equation*}
The unperturbed system is now Hamiltonian with a center at the
origin whose period annulus is bounded by a saddle loop. We have
written the transformations so as to directly  apply
\teoc{quadratic}. The Hamiltonian function of the unperturbed
system is
 \[
  H(x,y)=A(x)+B(x)y^2\,\mbox{ with $A(x)=\frac{1}{4}x^2(x+2)$ and $B(x)=x+2.$}
 \]
The projection of the period annulus is $(-4/3,2/3)$ and the
polycycle at its outer boundary has energy level $h_0=A(2/3)=8/27.$
It is very easy to show that the first Melnikov function is a linear
combination of
 \[
  \widetilde I_i(h)=\int_{\gamma_h}(x+2)^i y dx\,\mbox{ for $i=0,1,2.$}
 \]
Hence, the aforementioned result will follow once we check that
$\bigl\{\widetilde I_0,\widetilde I_1,\widetilde I_2\}$ is an
ECT-system on $(0,h_0).$  By using \lemc{puja} exactly as before,
$\widetilde I_i(h)=\frac{1}{h}I_i(h)$ where
 \[
 I_i(h)=\int_{\gamma_h}f_i(x)y^3dx
 \]
with $f_0(x)=\frac{2(x+2)(15{x}^{2}+42x+32)}{3(3x+4)^{2}},$
$f_1(x)=\frac{4(x+2)^2(9{x}^{2}+23x+16)}{3(3x+4)^{2}}$ and
$f_2(x)=\frac{2(x+2)^3(21{x}^{2}+50x+32)}{3(3x+4)^{2}}.$ Once
again, $\bigl\{\widetilde I_0,\widetilde I_1,\widetilde
I_2\bigr\}$ is an ECT-system on $(0,h_0)$ if, and only if, so it
is $\{I_0,I_1,I_2\}.$ The involution associated to~$A$ is
$z=\sigma(x)$ given by $q(x,z)\!:=x^2+xz+2x+z^2+2z=0$ because
$A(x)-A(z)=\frac{1}{4}(x-z)q(x,z).$ Thus
 \[
  \sigma'(x)=\frac{dz}{dx}=-\frac{z+2x+2}{x+2z+2}
 \]
and, setting $\ell_i=\mathscr
B_{\sigma}\!\left(\frac{f_i}{A'B^{3/2}}\right),$ we have to verify
that $W[\mathbf{\ell_i}]$ does not vanish on $(0,2/3)$ for
$i=1,2,3.$ It can be shown that, for $i=1,2,3,$
$W[\mathbf{\ell_i}](x)=\omega_i\bigl(x,\sigma(x)\bigr)$ with
$\omega_i(x,z)$ being a \emph{rational} function of $u=\sqrt{x+2}$
and $v=\sqrt{z+2},$ say $R_i(u,v).$ We note that
$x\longmapsto\sqrt{x+2}$ maps $(0,2/3)$ to $(\sqrt{2},\sqrt{8/3}).$
The resultant with respect to $v$ between the numerator of
$R_i(u,v)$ and $q(u^2-2,v^2-2)$ is a polynomial $r_i(u)$ that, by
applying Sturm's Theorem, has no roots on $(\sqrt{2},\sqrt{8/3}).$
Therefore, $W[\,\mathbf{\ell_i}]$ does not vanish on $(0,2/3)$ for
$i=1,2,3.$ By \teoc{quadratic}, we have proved the result of Peng in
\cite{Peng}.
\end{ex}

\subsection{Results on the program of Gautier, Gavrilov and
Iliev}\label{Gautier}

Our last examples of application come from the paper of Gautier,
Gavrilov and Iliev~\cite{Iliev3}, where a program for finding the
cyclicity of the period annuli of quadratic systems with centers of
genus one is presented. They give a list of the essential
perturbations of these centers (i.e., the one-parameter
perturbations that produce the maximal number of limit cycles),
together with the corresponding generating function of limit cycles
(i.e., the Poincar\'e-Pontryagin-Melnikov function). Since some
cases have been already solved in the literature about the problem,
this list includes only the open cases, a total of 26. They
conjecture that the cyclicity of these period annuli is two, except
for some particular cases in which it is three (cf. Conjecture~1 in
page~12 and Conjecture~2 in page~17). In their Theorem~3, two
quadratic reversible systems with a center are considered, denoted
by (r11) and (r18) in the list, and they show that, in both cases,
the upper bound of the number of limit cycles produced by the period
annulus under quadratic perturbations is equal to two. We are going
to reobtain this result for the case (r11) by using our criterion.
Moreover, we prove their conjecture in four new cases in their list,
namely (r7-r14), (r15), (r17) and (rlv3). In fact, Theorem
\ref{quadratic} is likely to be applied in many of their cases but
we have only been able to directly show that the functions on the
integrand satisfy the Chebyshev condition in the five mentioned
cases. We remark that our criterion gives a sufficient condition for
the Abelian integrals to be an ECT-system.

\medskip \noindent\textbf{Case (r11)} We translate the center to the
origin, so that the first integral of the unperturbed system is

 \[
  H(x,y)=A(x)+B(x)y^2\,\mbox{ with $A(x)=\frac{x^2(x+3)}{6(x+1)^3}$ and $B(x)=\frac{1}{2(x+1)^3}.$}
 \]
They show that the cyclicity of the period annulus under quadratic
perturbations is two. This will follow once we show that
$\bigl\{\widetilde I_0,\widetilde I_1,\widetilde I_2\bigr\}$ is an
ECT-system on $(0,1/6),$ where
 \[
  \widetilde I_i(h)=\int_{\gamma_h}(x+1)^{i-2}ydx.
 \]
The projection of the period annulus of the center at the origin is
$(-1/3,+\infty).$ By applying \lemc{puja} once again, $\widetilde
I_i(h)=\frac{1}{36 h}I_i(h)$ where
$I_i(h)=\int_{\gamma_h}f_i(x)y^3dx$ with
 \[
  f_0(x)=\frac{5{x}^{2}+13x+24}{(x+1)^{5}},\;
  f_1(x)=\frac{7{x}^{2}+19x+24}{(x+1)^{4}}\,\mbox{ and }
  f_2(x)=\frac{9{x}^{2}+25x+24}{(x+1)^{3}}.
 \]
It is clear then that it suffices to show that $\{I_0,I_1,I_2\}$
is an ECT-system on $(0,1/6).$ With this aim in view, let us note
that $A(x)-A(z)=\frac{(x-z)q(x,z)}{6(x+1)^3(z+1)^3}$ with
$q(x,z)\!:=3x^2z+x^2+10xz+3x+3xz^2+z^2+3z,$ so that the involution
$z=\sigma(x)$ associated to $A$ satisfies
$q\bigl(x,\sigma(x)\bigr)=0.$ Taking this into account, we get
that
 \[
  \sigma'(x)=\frac{dz}{dx}=-\frac{x(z+1)^4}{z(x+1)^4}.
 \]
As before we must compute the Wronskians $W[\,\mathbf{\ell_i}](x)$
for $i=1,2,3,$ where $\ell_i=\mathscr
B_{\sigma}\!\left(\frac{f_i}{A'B^{3/2}}\right),$ and then show that
they do not vanish for $x\in (0,+\infty).$ In this case
$W[\mathbf{\ell_i}](x)=\omega_i\bigl(x,\sigma(x)\bigr)$ with
$\omega_i(x,z)$ being a \emph{rational} function of $u=\sqrt{x+1}$
and $v=\sqrt{z+1},$ say $R_i(u,v).$ The resultant with respect to
$v$ between the numerator of $R_i(u,v)$ and $q(u^2-1,v^2-1)$ is a
polynomial $r_i(u).$ Since the mapping $x\longmapsto\sqrt{x+1}$
sends $(0,+\infty)$ to $(1,+\infty),$ the result will follow once we
show that these polynomials $r_i(u)$ do not vanish on $(1,+\infty)$.
This latter fact is deduced from the application of Sturm's
Theorem.\qedB

Let us mention that we have studied the case (r18) as well (the
other case that contemplates Theorem~3 in~\cite{Iliev3}), but it
seems that it cannot be solved by using the criterion given by our
\teoc{quadratic}. Of course, the success in the application of this
criterion depends on the particular problem studied, but we want to
stress that, when it works, it enables to extremely simplify the
solution. For instance, the proof of Theorem~3 takes eight pages of
highly nontrivial arguments. From now on, for the sake of brevity in
the exposition, we omit many of the explanations on the way to apply
our criterion since they are a verbatim repetition of the previous
examples.

\medskip\noindent\textbf{Cases (r7-r14) and (r15)} The first
integral is shared by the two cases and, after we translate the
center at the origin, it reads for
\[
 H(x,y)=\frac{y^2}{2}+\frac{x^2(3x^2+8x+6)}{12}.
\]
The cyclicity of the period annulus, whose projection on the
$x$-axis is the interval $(-1,1/3)$, is two if we prove that
$\bigl\{\widetilde I_0,\widetilde I_1,\widetilde I_2\bigr\}$ is an
ECT-system for $h\in (0,1/12),$ where
\[ \begin{array}{l}
\displaystyle \widetilde I_i(h)\, =\, \int_{\gamma_h} (x+1)^{i-2}
y dx \ \ \mbox{for the case (r7-r14)},
\vspace{0.2cm} \\
\displaystyle \widetilde I_i(h)\, =\, \int_{\gamma_h} (x+1)^{i-4}
y dx \ \ \mbox{for the case (r15)}.
\end{array}
\]
We apply \lemc{puja} to the Abelian integrals given by $ I_i(h)=h
\, \widetilde I_i(h)$ in order to write them in the form
$I_i(h)=\int_{\gamma_h}f_i(x) y^3 dx$. We have that:
\[
\begin{array}{l}
\displaystyle \int_{\gamma_h} H(x,y) y dx \, = \, \int_{\gamma_h}
\frac{21x^3+63x^2+64x+24}{36(x+1)^3} \, y^3 dx,
\\[10pt]
\displaystyle \int_{\gamma_h} H(x,y) \, \frac{ydx}{x+1}=
\int_{\gamma_h} \frac{(2x+3)(9x^2+14x+8)}{36(x+1)^4} \, y^3 dx,
\\[10pt]
\displaystyle \int_{\gamma_h} H(x,y) \, \frac{ydx}{(x+1)^2}=
\int_{\gamma_h} \frac{15x^3+47x^2+52x+24}{36(x+1)^5} \, y^3 dx,
\\[10pt]
\displaystyle \int_{\gamma_h} H(x,y) \, \frac{ydx}{(x+1)^3}=
\int_{\gamma_h} \frac{12x^3+39x^2+46x+24}{36(x+1)^6} \, y^3 dx,
\\[10pt]
\displaystyle \int_{\gamma_h} H(x,y) \,
\frac{ydx}{(x+1)^4}=\int_{\gamma_h}
\frac{9x^3+31x^2+40x+24}{36(x+1)^7} \, y^3 dx.
\end{array}
\]
Some computations show that the involution $\sigma$ defined by
$A(x)\!:=H(x,0)$ satisfies $q\bigl(x,\sigma(x)\bigr)=0$ with
$q(x,z)\!:= 3z^3+3xz^2+8z^2+3x^2z+8xz+6z+3x^3+8x^2+6x$. We use
resultants and Sturm's Theorem in order to check that the
corresponding Wronskians have no zeros on the interval
$(0,1/3)$.\qedB

\medskip\noindent\textbf{Case (r17)} Once the center is translated to the
origin, the first integral reads for
\[
 H(x,y)=\frac{y^2}{2}+\frac{(2x+3)x^2}{6}.
\]
Setting $\widetilde I_i(h)=\int_{\gamma_h} (x+1)^{i-3} y dx,$ the
cyclicity of its period annulus is two if we prove that
$\bigl\{\widetilde I_0,\widetilde I_1,\widetilde I_2\bigr\}$ is an
ECT-system on $(0,1/6).$ By \lemc{puja}, we have that $\widetilde
I_i(h)=\frac{1}{18h}\int_{\gamma_h}f_i(x) y^3 dx,$ with
\[
f_0(x)=\frac{5x^2+13x+12}{(x+1)^5}, \
f_1(x)=\frac{7x^2+16x+12}{(x+1)^4}\ \mbox{and} \ f_2(x)=
\frac{9x^2+19x+12}{(x+1)^3}.
\]
In this case, the involution $\sigma$ defined by $A(x)\!:=H(x,0)$
satisfies $q\bigl(x,\sigma(x)\bigr)=0$ where $q(x,z)\!:=
2z^2+2xz+3z+2x^2+3x$. The projection of the period annulus on the
$x$-axis is $(-1,1/2)$ and, thus, we are done if we show that the
functions $\ell_i=\mathscr
B_{\sigma}\!\left(\frac{f_i}{A'B^{3/2}}\right)$ form an ECT-system
in $(0,1/2)$. Once again, the involution can be explicitly written,
but we prefer to use resultants and Sturm's Theorem because it
provides an algebraic procedure to check that the Wronskians
$W[\,\mathbf{\ell_i}]$ do not vanish on $(0,1/2)$ for $i=1,2,3.$ The
proof of this fact is omitted for the sake of shortness.\qedB

\medskip\noindent\textbf{Case (rlv3)} After the center is translated to the
origin, the first integral becomes
 \[
  H(x,y)=x^2(2-x^2)+\frac{1}{2}(1+x)^2y^2.
 \]
Since $A(x)\!:=H(x,0)$ is an even function, we have that
$\sigma(x)=-x$ and this simplifies a lot the computations. The
projection of the period annulus on the $x$-axis is $(-1,1)$. In
order to prove that its cyclicity under quadratic perturbations is
two, we are lead to show that $\bigl\{ I_0, I_1, I_2\bigr\}$ form
an ECT-system for $h\in (0,1),$ where $I_i(h)=\int_{\gamma_h}
f_i(x) y^3 dx$ with $f_0(x)=
\frac{(5x^4-2x^3-9x^2+4x+8)(x+1)}{2(x-1)^4},$ $f_1(x)=
\frac{(7x^4-13x^2+8)(x+1)}{(x-1)^2}$ and $f_2(x)=
\frac{6x^4+x^3-11x^2-2x+8}{(x-1)^2}.$ To this end, by applying
\teoc{quadratic} and taking $\sigma(x)=-x$ into account, it
suffices to show that the functions
\[
 \ell_0(x)=\frac{5x^6- 8x^4+ 7 x^2+8}{2x(x-1)^5(x+1)^5}, \ \
 \ell_1(x)=\frac{ 7 x^4- 13 x^2+8}{x(x-1)^3(x+1)^3} \ \mbox{ and } \
 \ell_2(x)=\frac{5 x^4- 9 x^2+8}{x(x-1)^4(x+1)^4}
\]
form an ECT-system on $(0,1).$ It is easy to see that $\ell_2$ does
not vanish on $(0,1)$. The Wronskian associated to $\ell_1$ and
$\ell_2$ is the rational function
\[
W[\ell_1,\ell_2](x)=\frac{96 - 240 x^2 + 243 x^4 - 126 x^6 + 35
x^8}{18x(x-1)^8(x+1)^8},
\] which has no zero on $(0,1)$ by
virtue of Sturm's Theorem. Finally
\[
 W[\ell_0,\ell_1,\ell_2](x)=\frac{2(512 - 1632x^2 + 2200 x^4 - 1617 x^6 + 693 x^8 -
175 x^{10} + 35 x^{12})}{9(x-1)^{15}(x+1)^{15}},
\]
which neither vanishes on $(0,1)$, again by using Sturm's Theorem.
As desired, this shows that $(\ell_2,\ell_1,\ell_0)$ is an
ECT-system on $(0,1).$\qedB

\section{Appendix}\label{Sec5}

\subsection{Resultant of two polynomials}

Given two polynomials $p,q\in \mathbb{C}[x,y],$  say
 \begin{align*}
  &p(x)=a_0x^m+a_1x^{m-1}+\ldots+a_m,\,\mbox{ with $a_0\neq 0,$} \\
  &q(x)=b_0x^n+b_1x^{n-1}+\ldots+b_n,\,\mbox{ with $b_0\neq 0,$}
 \end{align*}
where $a_i, b_i \in \mathbb{C}[y]$, the \emph{resultant} of $p$
and $q$ with respect to $x$, denoted by $\mbox{Res}(p,q,x)$ is the
$(m+n)\!\times\!(m+n)$ determinant
 \[
  \mbox{Res}(p,q,x) \, =\, \det
\left(
  \begin{array}{cccccccc}
    a_0 &  &  &  & b_0 &  &  &  \\
    a_1 & a_0 &  &  & b_1 & b_0 &  &  \\
    a_2 & a_1 & \ddots &  & b_2 & b_1 & \ddots &  \\
    \vdots & a_2 & \ddots & a_0 & \vdots & \ddots &  & b_0 \\
    a_m & \vdots & \ddots & a_1 & b_n & \vdots & \ddots & b_1 \\
     & a_m &  & a_2 &  & b_n &  & b_2 \\
     &  & \ddots & \vdots &  &  & \ddots & \vdots \\
     &  &  & a_m &  &  &  & b_n \\
  \end{array}
\right)
 \]
where the blank spaces are filled with zeros. The three basic
properties of the resultant are:
\begin{enumerate}
\item $\mbox{Res}(p,q,x)$ is an integer polynomial in the
      coefficients of $p$ and $q$. \item $\mbox{Res}(p,q,x)=0$ if, and
      only if, $p$ and $q$ have a nontrivial common factor in $\mathbb{C}[x,y]$.
\item There are polynomials $A,B\in \mathbb{C}[x,y]$ such that
      $A p + B q = \mbox{Res}(p,q,x).$ Moreover
      the coefficients of~$A$ and~$B$ are integer polynomials in the coefficients of $p$
      and $q$.
\end{enumerate}

Resultants can be used to eliminate variables from systems of
polynomial equations. As an example, let us suppose that we want
to study the following system of two polynomial equations with two
variables:
 \[
  \left\{
  \begin{array}{l}
   xy-1=0, \\[2pt]
   x^2+y^2-4=0.
  \end{array}
  \right.
 \]
Here we have two variables to work with, but if we regard
$p(x,y)\!:=xy-1$ and $q(x,y)\!:=x^2+y^2-4$ as polynomials in~$x$
whose coefficients are polynomials in~$y,$ we can compute the
resultant with respect to $x$ to obtain
$\mbox{Res}(p,q,x)=y^4-4y^2+1.$ By the third property above, there
are polynomials $A,B\in\C[x,y]$ such that
$A(x,y)p(x,y)+B(x,y)q(x,y)=y^4-4y^2+1.$ Accordingly, $y^4-4y^2+1$
vanishes at any common solution of $p=q=0.$ Thus, we can solve
$y^4-4y^2+1=0$ and find the $y$-coordinates of these solutions.

\subsection{Sturm's Theorem}

A sequence $\{f_0,f_1,\ldots,f_{m}\}$ of continuous real functions
on $[a,b]$ is called a \emph{Sturm's sequence} for $f=f_0$ on
$[a,b]$ if the following is verified:
 \begin{enumerate}
  \item $f_0$ is differentiable on $[a,b].$
  \item $f_m$ does not vanish on $[a,b].$
  \item If $f(x_0)=0$ with $x_0\in [a,b]$ then $f_1(x_0)f_0'(x_0)>0.$
  \item If $f_i(x_0)=0$ with $x_0\in [a,b]$ then $f_{i+1}(x_0)f_{i-1}(x_0)<0.$
 \end{enumerate}

\begin{proclama}{Sturm's Theorem.}
Let $\{f_0,f_1,\ldots,f_{m}\}$ be a Sturm's sequence for $f=f_0$ on $[a,b]$ with
$f(a)f(b)\neq 0.$ Then the number of roots of $f$ on $(a,b)$ is equal to $V(a)-V(b),$
where $V(c)$ is the number of changes of sign in the sequence
$\{f_0(c),f_1(c),\ldots,f_m(c)\}.$
\end{proclama}

There is a simple procedure to construct a Sturm's sequence in
case that $f$ is polynomial. Indeed, if $p(x)$ is a polynomial of
degree $n$, we define the sequence $\{p_0,p_1,\ldots,p_m\}$ with
$m\leqslant n$ in the following way. We set $p_0=p,$ $p_1=p'$ and
 \begin{align*}
  & p_{i-1}(x)=q_i(x)p_i(x)-p_{i+1}(x),\,\mbox{ for $i=1,2,\ldots,m-1,$} \\[2pt]
  & p_{m-1}(x)=q_m(x)p_m(x),
 \end{align*}
where $q_i(x)$ and $p_{i+1}(x)$ are the quotient and the remainder
(the latter with the sign changed) of the division of $p_{i-1}(x)$
by $p_{i}(x)$, respectively. The construction of this sequence
ends when the remainder is zero, i.e., $p_{m+1}=0.$ In this case,
since this is essentially Euclides' algorithm, $p_m$ is the
greatest common divisor of $p_0$ and $p_1.$ If all the zeros of
$p$ are simple then $p_m$ does not vanish and it is easy to show
that $\{p_0,p_1,\ldots,p_m\}$ is a Sturm's sequence for $p$ on any
interval. If $p$ has zeros with multiplicity then $p_m$ vanishes.
Since $p_m$ divides $p_0$ and $p_1,$ it also divides $p_i$ for
$i=2,3,\ldots,m$. In this case, we set $\bar p_i=p_i/p_m$ and it
follows that $\{\bar p_0,\bar p_1,\ldots,\bar p_m\}$ is a Sturm's
sequence for $p$ on any interval.

\bibliographystyle{plain}

\end{document}